\documentclass[11pt]{llncs}
\usepackage{graphicx}

\usepackage{times}
\usepackage{mathptm}

\usepackage[round,sort]{natbib}
\DeclareSymbolFont{AMSb}{U}{msb}{m}{n}
\DeclareSymbolFontAlphabet{\Bbb}{AMSb}
\def\R{\ensuremath{\Bbb R}}

\makeatletter
\def\hb@xt@{\hbox to }
\makeatother

\let\oldendproof\endproof
\def\endproof{\qed\oldendproof}

\begin{document}

\title{Manhattan Orbifolds} 

\author{David Eppstein}

\institute{Computer Science Department\\
University of California, Irvine\\
\email{eppstein@uci.edu}}

\maketitle   

\begin{abstract}
We investigate a class of metrics for 2-manifolds in which, except for a discrete set of singular points, the metric is locally isometric to an $L_1$ (or equivalently $L_\infty$) metric, and show that with certain additional conditions such metrics are injective. We use this construction to find the tight span of squaregraphs and related graphs, and we find an injective metric that approximates the distances in the hyperbolic plane analogously to the way the rectilinear metrics approximate the Euclidean distance.
\end{abstract}

\section{Introduction}

An \emph{injective metric space} is a metric space $X$ such that whenever $X$ is isometric to a subset of some larger metric space $Y$, there exists a nonexpansive mapping from $Y$ to itself that fixes that subset and maps the rest of $Y$ onto it. Equivalently by a theorem of \citet{AroPan-PJM-56} an injective space is a path-geodesic metric space in which any family of closed metric balls forms a Helly family: for any clique in an intersection graph of balls, there is a point in the space contained in all the balls of the clique. Standard examples of injective spaces include the $L_{\infty}$ metric on any real vector space, the $L_1$ or \emph{Manhattan metric} on the plane (equivalent to the $L_{\infty}$ metric by rotation and scaling), and the metric on any real tree.

As \citet{Isb-CMH-64} showed, any metric space can be isometrically embedded in a unique minimal injective space called its \emph{injective envelope}, \emph{hyperconvex hull}, or \emph{tight span}. Tight spans have gained attention recently for their applications in the reconstruction of evolutionary trees \citep{DreHubMou-DM-01} and in online algorithms \citep{ChrLar-Algs-94,BeiChrLar-TCS-02,cs.DS/0611088}.

In this paper we investigate injective metric spaces that have the topology of a 2-manifold. In particular, we describe a class of such spaces, which we call {\em Manhattan orbifolds} as they are modeled after the Manhattan metric except at a discrete set of singular points similar to the elliptic singularities of Thurston's orbifolds.

We then use these surfaces to construct the tight spans of certain planar graphs, the \emph{squaregraphs} and some related nonbipartite and nonplanar graphs. A squaregraph is a planar graph in which all but one face in some embedding are quadrilaterals, and any vertex that is not on the non-quadrilateral face has degree at least four. As we show, filling each quadrilateral face of a squaregraph with a unit square in the $L_1$ metric produces a Manhattan orbifold, the so-called $L_1$ \emph{median complex} of the graph, which is the tight span of the graph.
We describe the result of applying this construction to an infinite squaregraph, the $\{4,.5\}$ tesselation of the hyperbolic plane. The result is an injective space with the same topology as the hyperbolic plane, in which the hyperbolic distance between any two points is within a constant factor of the distance between the same two points in the injective metric.

Finally, we close with an application of Manhattan orbifolds in describing \emph{greedy embeddings}~\cite{PapRat-TCS-05} of arbitrary graphs.

\section{Rectilinear cones}

\begin{figure}[t]
\centering\includegraphics[height=3in]{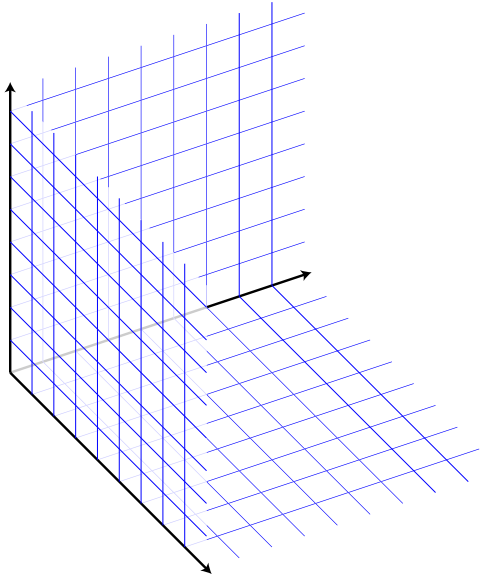}
\caption{The order-3 rectilinear cone, formed by gluing together three quadrants of a Manhattan metric plane.}
\label{fig:order-3-cone}
\end{figure}

The spaces we study in this paper are locally modeled after the \emph{Manhattan metric} in the plane $\R^2$, that is, the $L_1$ metric in which the distance between any two points $(x_1,y_1)$ and $(x_2,y_2)$ is $|x_1-x_2|+|y_1-y_2|$. This metric is also equivalent (by rotation and scaling) to the \emph{Chebyshev metric} (the $L_\infty$ metric) in which the distance is instead $\max(|x_1-x_2|,|y_1-y_2|)$. We say that a space is \emph{locally Manhattan} at a point $p$ if some neighborhood of $p$ is isometric to a neighborhood of a point in the Manhattan plane. We begin by describing the spaces that are locally Manhattan except at a single point of singularity, a \emph{cone point}.

Consider the (non-injective) two-dimensional surface formed by the boundary of the positive orthant in $\R^3$ (Figure~\ref{fig:order-3-cone}), with the metric on the surface being the $L_1$ metric for $\R^3$. Any two points in this surface can be connected by a path composed of axis-aligned line segments, lying in the surface, with length equal to the distance between the two points, so this surface is path-geodesic. Geometrically, this surface is composed of three infinite right-angled plane wedges, each isometric to the positive quadrant of the $L_1$ plane, glued together along their boundary rays. However, in the intrinsic geometry of the surface, there is nothing special about the points along which pairs of quadrants are glued: in a neighborhood of any such point, the metric is the same as it would be in a neighborhood of a point on one of the coordinate axes of the $L_1$ plane, where its four quadrants are glued together; but these points are not different than any other point of the $L_1$ plane. That is, this surface is locally Manhattan everywhere except at the origin of $\R^3$, where there are three quadrants meeting while everywhere else there are four. We say that the origin is a \emph{cone point} of this surface.

Similarly, for any $k>4$ we can form a surface, the \emph{order-$k$ rectilinear cone}, by gluing together $k$ quadrants of the Manhattan plane, along their boundary rays, and letting the distance between any two points of this surface equal the length of the shortest path connecting them. Such a surface can be embedded isometrically into the $L_1$ metric for $\R^k$, with the gluing rays on orthogonal coordinate axes. It is locally Manhattan except at a cone point, the point forming the common origin of the glued quadrants. We define the \emph{angular excess} of the cone point to be $2\pi-k(\pi/2)$.

More generally, we define a \emph{cone point} of a metric space to be any point that has a neighborhood isometric to a neighborhood of the origin in an order-$k$ rectilinear cone, and we define the angular excess of such a point to be $2\pi-k(\pi/2)$. Any cone point has a unique $k$ satisfying this definition; this can be seen, for instance, from the fact that a sufficiently small metric ball around a cone point has a boundary in the form of a polygon with $k$ right angles.

\section{Boundary singularities}

\begin{figure}[t]
\centering\includegraphics[height=1.75in]{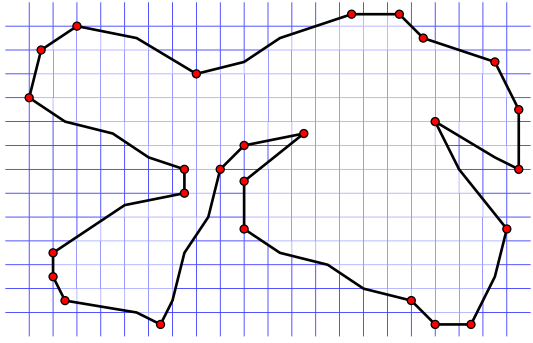}
\caption{A region of the Manhattan plane with piecewise smooth boundary, showing the inflection points on the boundary.}
\label{fig:boundary-geodesy}
\end{figure}

Along with cone points in the interior of a manifold, we also need to model certain kinds of singularities on the boundary of the manifold. We take as our model a region $R$ of the Manhattan plane with piecewise smooth boundary, such as a simple polygon (Figure~\ref{fig:boundary-geodesy}). We partition the smooth boundary points into eight subsets: those with slopes $0$, $1$, $-1$, and $\infty$, and those with slopes in the four open intervals $(\infty,-1)$, $(-1,0)$, $(0,1)$, and $(1,\infty)$ between these four slopes. We say that a boundary point $b$ of such a region is \emph{boundary-geodesic} if there is a neighborhood $N$ of $b$ such that, within $N$, all boundary points of $R$ are smooth and all belong to a single one of these eight subsets. If $b$ does not have this property, but it has a neighborhood within which all boundary points of $R$ except $b$ itself are boundary-geodesic, we call $b$ an {\em inflection point}. Note that a connected set of boundary-geodesic points forms a geodesic in the Manhattan plane. Note also that these definitions do not depend on the choice of axes for the plane, as rotating the plane by a multiple of $\pi/2$ or reflecting it across one of the coordinate axes will preserve the partition of boundary points according to their slopes. In the figure, the inflection points are marked as small red circles, and the remaining boundary points are boundary-geodesic.

We may extend these definitions to metric spaces locally modeled on regions of the Manhattan plane: We define a point $b$ of a metric space to be boundary-geodesic if there exists an open neighborhood $N$ of $b$, a region $R$ bounded by a piecewise smooth curve in the Manhattan plane, and an isometry that maps $N$ to an open subset of $R$ and that maps $b$ to a boundary-geodesic point of $R$. Similarly, we define a point $b$ of a metric space to be an inflection point if there exists an open neighborhood $N$ of $b$, a region $R$ bounded by a piecewise smooth curve in the Manhattan plane, and an isometry that maps $N$ to an open subset of $R$ and that maps $b$ to an inflection point of $R$.

For the generalization to 2-manifolds, we will also need another type of singular point. We say that a point $b$ of a metric space is a \emph{cone inflection point} if there exists an open neighborhood $N$ of $b$, an integer $k$, a region $R$ bounded by a piecewise-smooth curve in an order-$k$ rectilinear cone, and an isometry that maps $N$ to an open subset $N'$ of $R$, such that the isometry maps $b$ to the origin of the cone and such that all points in $N'$ other than the origin are either locally Manhattan or boundary-geodesic. If $b$ is a cone inflection point for an order-$k$ rectilinear cone, it is also such a point for any rectilinear cone of order greater than $k$; this ambiguity will not cause us any difficulty.

\section{Manhattan orbifolds}

We are now ready to define Manhattan orbifolds, the surfaces that we will later show to be injective.
We define a Manhattan orbifold to be a Cauchy-complete metric space with the following properties:
\begin{itemize}
\item Every point has a neighborhood homeomorphic to an open disk or to an open semidisk. That is, the space is a 2-manifold with boundary.
\item Every simple closed curve is the boundary of a unique disk. That is, the surface obeys the \emph{Jordan curve theorem}: it can have no uncontractible cycles, and it cannot be homeomorphic to a sphere (because in a sphere, a simple closed curve bounds two disks). We call the disk bounded by a simple closed curve the \emph{interior} of the curve.
\item Every point with neighborhood homeomorphic to an open disk is either locally Manhattan, or a cone point with negative angular excess.
\item Every point with neighborhood homeomorphic to an open semidisk is boundary-geodesic, an inflection point, or a cone inflection point.
\item The distance between any two points in the space equals the length of the shortest curve connecting the two points, as measured in the local neighborhoods described above. That is, the space is path-geodesic.
\end{itemize}

Examples of Manhattan orbifolds include the Manhattan plane itself, any rectilinear cone of order $k$ for $k>4$, and any polygonal subset of these spaces. We will later describe some more complicated examples.

It may be of interest to consider spaces defined in a similar way that contain cone points of order three, or that contain uncontractible cycles; however such spaces cannot be injective and we do not further investigate them here.

\section{Singularities in a bounded region}

We define a \emph{singular point} of a Manhattan orbifold to be a point that is neither locally Manhattan nor boundary-geodesic. We observe that the singular points of a Manhattan orbifold must form a discrete subset: every point $p$ of the orbifold must have a neighborhood containing no singular point other than possibly $p$ itself, because all of the allowed types of points in the orbifold are defined in terms of neighborhoods that have no other singular point.

Our first technical result uses K\"onig's lemma \citep{Koe-36} to limit the number of singular points that may occur within any bounded region of a Manhattan orbifold.

\begin{figure}[t]
\centering\includegraphics[height=2in]{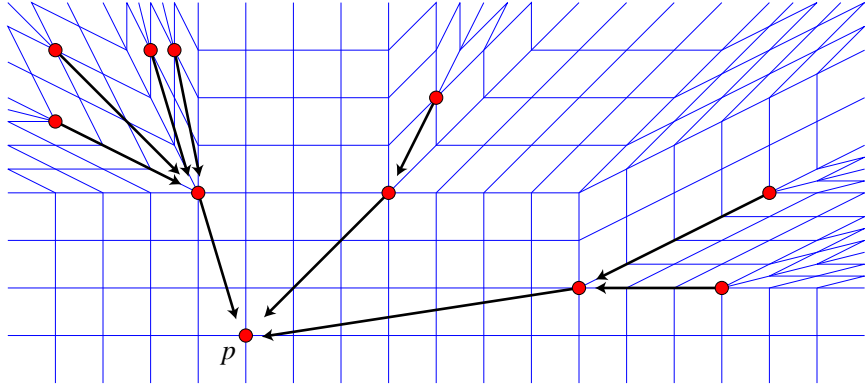}
\caption{The tree defined in the proof of Lemma~\ref{lem:finite-singular}, in which each singular point is connected to a nearer singular point to $p$, or to $p$ itself if no such point exists.}
\label{fig:koenig}
\end{figure}

\begin{lemma}
\label{lem:finite-singular}
For any point $p$ of a Manhattan orbifold $M$, and any bounded radius $r>0$, finitely many points within distance $r$ of $p$ are singular.
\end{lemma}

\begin{proof}
Let $S$ denote the set of singular points within distance $r$ of $p$, together with $p$ itself. Form a directed graph $G$ having the points of $S$ as vertices. For each point $s$ in $S\setminus\{p\}$, add to $G$ a single edge $(s,t)$, where $t$ is chosen to have the properties that there is a geodesic from $s$ to $p$ that passes through $t$ but that there does not exist a geodesic from $s$ to $t$ that passes through any other point in $S$. Such a point $t$ may be found by letting $t_0=p$ and then (for $i>0$) letting $t_i$ be any point of $S$ that lies on a geodesic from $s$ to $t_{i-1}$, until the sequence can not be extended any longer. The points in this sequence $t_i$ have decreasing distances to $s$, so the sequence must eventually terminate, for if there were infinitely many points in the sequence then (by Cauchy closure) there would be a limit point of the sequence with infinitely many singular points in its neighborhoods, violating the definition of a Manhattan orbifold. The final point in the sequence may be chosen as $t$. If multiple alternative choices for $t$ exist, choose one arbitrarily. Figure~\ref{fig:koenig} depicts an example of this construction; in the figure, the lines parallel to axes of the locally defined Manhattan metrics are indicated by the light blue lines, and the singular points (shown as red circles) are the ones with more than four such axes incident to them.

Then, for any $s\ne p$ in $S$, there is a unique path of outgoing edges in $G$ that begins at $s$. Each step on this path reduces the remaining distance to $p$, so the total geodesic length of such a path is at most $r$. By the Cauchy-completeness of $M$ and the discreteness of $S$, such a path cannot continue for more than finitely many steps, so it must eventually reach the only vertex without an outgoing edge, $p$. Therefore, $G$ is a tree, rooted at $p$. Any path in $G$ has geodesic length equal to the distance between the path's endpoints, at most $r$. Therefore, the graph formed by reversing all edges in $G$ has no path outwards from $p$ that contains infinitely many edges and vertices, because (again invoking Cauchy-completeness and discreteness) the vertices of such a path would have to have an accumulation point violating the requirement that all points of a Manhattan orbifold have a neighborhood without other singularities than the point itself.

We now examine the cardinality of the set of incoming edges of $G$ that may exist at any vertex $s$. Any neighborhood of $S$ can be partitioned into $k$ regions isometric to a neighborhood of the origin in a single quadrant of the Manhattan plane, where $k$ is the order of $s$ as a cone point or cone inflection point (or $4$, if $s$ is not a point of that type). We partition the incoming neighbors of $s$ according to which of these quadrants contain the geodesic from the neighbor to $s$, and only consider points in one set of this partition at a time.

For any edge $(t,s)$ in $G$, the set of points in $M$ that belong to geodesics from $t$ to $s$ is isometric to a rectangle of the Manhattan plane; it can have no singularities within it because if such a singularity existed $t$ would be connected to such a singularity instead of to $s$. Thus, the union of these rectangles for all neighbors $t$ in a single quadrant of $s$ forms a subset of $M$ isometric to a union of rectangles in the Manhattan plane. The neighbors of $s$, and $s$ itself, lie on the boundary of this subset, and each point of the subset is within distance $r-d(s,p)$ of $s$. Any discrete bounded subset of the Manhattan plane is finite, so the set of neighbors of $s$ within a single quadrant is finite.

Since $s$ has finitely many quadrants, each containing finitely many neighbors, it has finitely many neighbors overall. Since $G$ is a tree with no infinite path and with finitely many children per vertex, by K\"onig's lemma, it is itself finite.
\end{proof}

\section{Orthogonal polygons}

\begin{figure}[t]
\centering\includegraphics[height=1.5in]{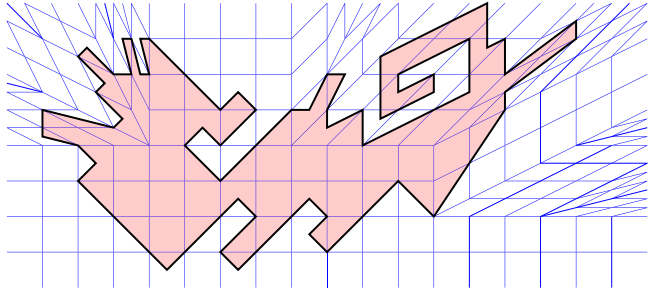}
\caption{An orthogonal polygon.}
\label{fig:orthog-poly}
\end{figure}

We define an \emph{orthogonal polygon} in a Manhattan orbifold $M$ to be a simple closed curve in $M$, of bounded total length, such that all but finitely many points of the curve have a neighborhood in $M$
that can be mapped isometrically to the Manhattan plane in such a way that the part of the curve within the neighborhood is mapped to a line of slope 1 in the Manhattan plane; see Figure~\ref{fig:orthog-poly} for an example. We call the finitely many remaining points at which this property does not hold the \emph{vertices} of the polygon. Recall that in a Manhattan orbifold any simple closed curve is required to bound a unique disk, its interior. At any vertex, we may define an \emph{interior angle}, an integer multiple $k$ of $\pi/2$, such that the portion of the interior of the polygon in a neighborhood of the vertex is isometric to $k$ quadrants of the $L_\infty$-metric plane glued together (or equivalently $2k$ half-quadrants of the Manhattan plane glued together). We may define the exterior angle similarly; note that, at a cone point, the interior and exterior angles do not add to $2\pi$.

\begin{lemma}
\label{lem:four-corners}
Any orthogonal polygon has at least four vertices with interior angle $\pi/2$.
\end{lemma}

\begin{proof}
The interior of the polygon has bounded radius, for if not the polygon would form an uncontractible  curve, contradicting the assumption that $M$ is a Manhattan orbifold and hence has no uncontractible closed curves. By Lemma~\ref{lem:finite-singular}, the interior of the polygon contains a finite number of singular points.

Let $S$ be the topological space formed by gluing together two copies of the interior of the polygon. We form a metric space from $S$ by replacing the Manhattan metric by the Euclidean metric for the same set of points, within any open region of $S$ that does not containing a singularity. Then $S$ is topologically a sphere, with locally Euclidean metric except at the singularities and boundary vertices, which have the same angular defects as they do in $M$. By the Gauss-Bonnet formula, the total angular defect of $S$ is $4\pi$. The only possible positive angular defect in $S$ is $\pi$, at a vertex of the polygon with interior angle $\pi/2$, so there must be at least four such vertices.
\end{proof}

\section{Orthoconvexity}

We wish to eventually prove the Helly property for balls in Manhattan orbifolds. These balls will (it turns out) resemble orthogonal polygons, but differ from them in two ways. First, orthogonal polygons may have large interior angles (concavities) at some of their vertices, a feature that is not possible for balls. And second, orthogonal polygons cannot include portions of the boundary of the orbifold that are not themselves isometric to slope-1 lines, while this restriction does not exist for balls. Therefore, it is convenient to define a class of shapes related to orthogonal polygons, but with restricted interior angles and with less restriction on how these shapes may meet the boundary of the orbifold. We define an \emph{orthoconvex region} to be a simply-connected bounded subset of a Manhattan orbifold, such that
\begin{itemize}
\item the boundary of the region meets the boundary of the orbifold in a finite number of components,
\item all but finitely many points of the boundary of the region either belong to the boundary of the orbifold or have a neighborhood in $M$ within which the curve is homeomorphic to a line of slope 1 in the Manhattan plane, and
\item at each of the remaining points of the boundary of the region, the interior angle is either $\pi/2$ or $\pi$.
\end{itemize}

\begin{figure}[t]
\centering\includegraphics[height=2in]{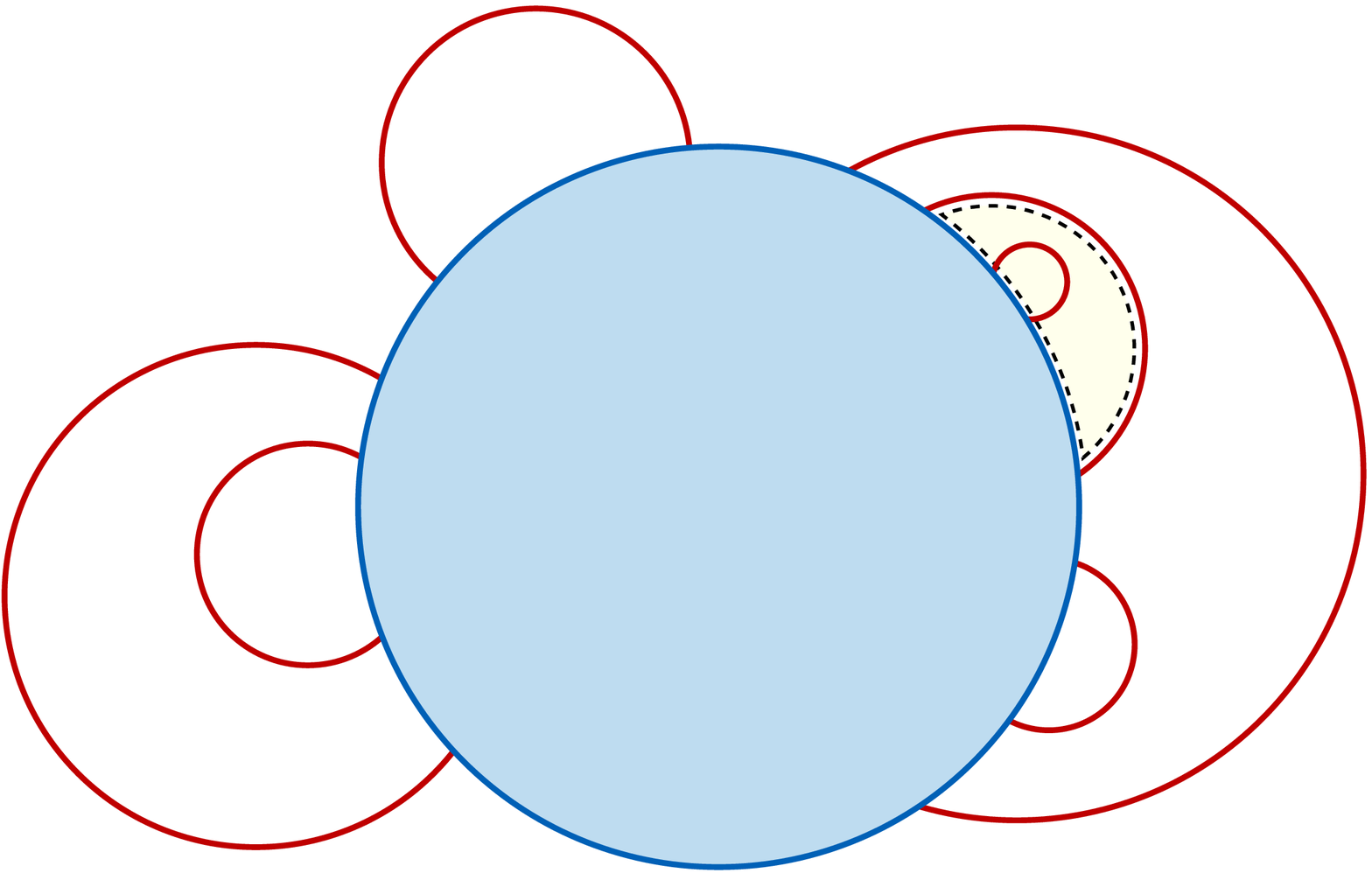}
\qquad
\includegraphics[height=2in]{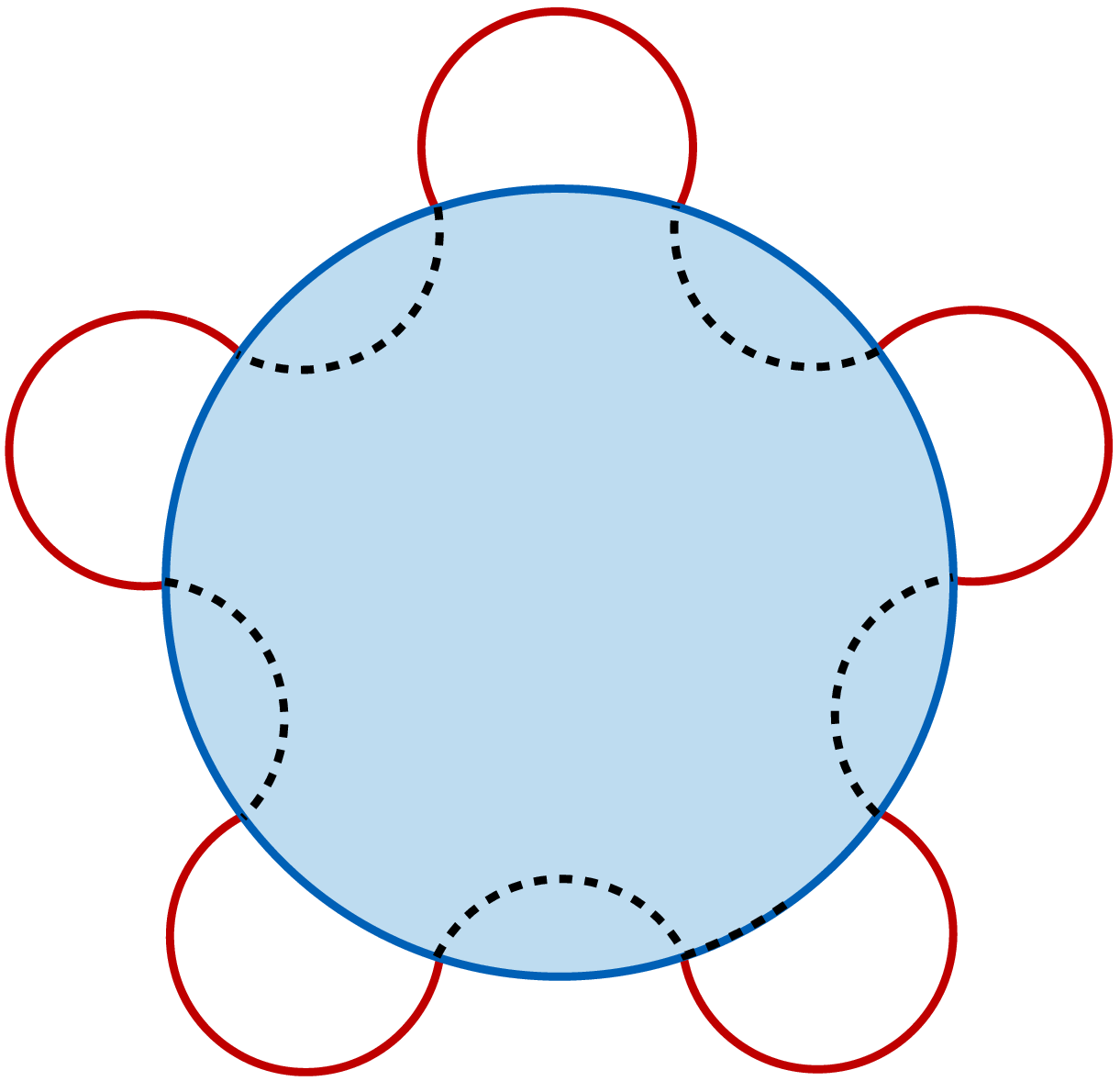}
\caption{Illustrations for the proof of Lemma~\ref{lem:ocintcon}, showing the boundary of the two orthoconvex regions $A$ and $B$ schematically, as topological curves. In both figures $A$ is the shaded region, and the other arcs denote portions of the boundary of $B$. Left: If two arc of the boundary of $B$ outside $A$ are nested, the inner one together with the boundary of $B$ forms an orthogonal polygon with too few right angles. Right: outer arcs must be connected in consecutive order to form the boundary of $B$.}
\label{fig:nesting-digon}
\end{figure}

In the remainder of this section we describe the intersection properties of orthoconvex regions.

\begin{lemma}
\label{lem:ocintcon}
Any nonempty intersection of any two orthoconvex regions is an orthoconvex region.
\end{lemma}

\begin{proof}
Let the two regions be $A$ and $B$. Since they are both topological disks, they have simple closed curves (which may include portions of the boundary of the manifold) as their boundaries. Define a \emph{arc} of $B$ to be a maximal subset of the boundary of $B$ that does not intersect $A$. Topologically, the arcs of $B$ form open curves that begin and end on the boundary of $A$ and are otherwise disjoint from $A$ (Figure~\ref{fig:nesting-digon}). The endpoints of any arc partition the boundary of $A$ in two parts; either of those two parts, together with the arc itself, forms a simple closed curve. It is not possible for two arcs with endpoints $c,d$ and $e,f$ to have those endpoints lie in the cyclic order $c,e,d,f$ around the boundary of $A$, because then only one of $e$ or $f$ would be interior to the closed curve formed by arc $cd$; therefore, any two arcs must form intervals on the boundary of $A$ that are either nested or disjoint, as shown in the left of the figure.

Define an arc $cd$ to be \emph{inner} if some other arc $ef$ forms (with $A$) a simple cycle such that $cd$ is interior to the cycle, and \emph{outer} otherwise. The outer arcs, together with the portions of the boundary of $A$ that are not interior to any arc's cycle, themselves form one large simple cycle with all of $A$ interior to it. $B$, too, must be interior to this cycle of outer arcs, because we have assumed that $A$ and $B$ intersect. Next, suppose there exists at least one inner arc, and let $cd$ be an inner arc that is contained in the cycle formed by an outer arc but that is not contained in the cycle formed by any other inner arc. Points $c$ and $d$ partition the boundary of $A$ into two parts; form a simple cycle from one of these parts and arc $cd$, choosing the part such that this simple cycle has $A$ exterior to it (the dashed cycle in the figure). Then $B$ must be exterior to this cycle as well, so its interior is a disk disjoint from $A\cup B$. In particular, this cycle cannot lie along any portion of the boundary of the manifold, because every point on the cycle has this inner disk on one side of it and $A$ or $B$ on the other side. Thus, it forms an orthogonal polygon. However, it can only have two internal angles of $\pi/2$, at $c$ or $d$, violating Lemma~\ref{lem:four-corners}. This contradiction shows that there can be no inner arcs.

Finally, suppose there are only outer arcs; that is, the arcs of $B$ enclose a cyclic sequence of disjoint intervals of the boundary of $A$. It remains to describe how these arcs may be connected to each other to form the boundary of $B$. The only possible connection pattern is that the right endpoint of each one of these arcs is connected to the left endpoint of the next arc, consecutively around the sequence of the outer arcs, as shown in the right of the figure. For, if any other two arc endpoints were to be connected, they would either separate an odd number of endpoints on each side of the connection, preventing the formation of a set of disjoint curves within $A$ connecting all the endpoints, or they would separate one subset of the arcs from another subset, violating the assumption that $A$ has a single boundary curve.

The intersection $A\cap B$ is then bounded by the connecting curves within $A$, together with the portions of the boundary of $A$ that are enclosed by each of the outer arcs of $B$. Since we have established that these curves form a single cyclic chain, the intersection $A\cap B$ has only one boundary component and can therefore have only one connected component. It is straightforward to verify that the boundary of this intersection satisfies the other requirements of an orthoconvex region.
\end{proof}

\begin{figure}[t]
\centering\includegraphics[height=2in]{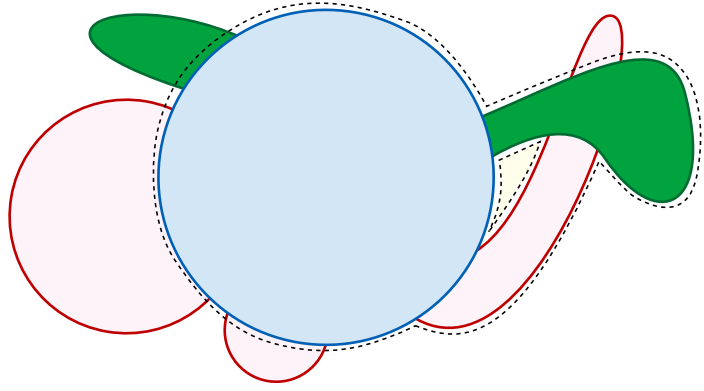}
\caption{Illustration for the proof of Lemma~\ref{lem:3way-helly}, showing schematically the three sets $A$ (medium blue shaded circle), $B$ (dark green), and $C$ (light red), together with the two closed curves formed by following part of an arc of $C$ and continuing around the boundaries of $B$ and $A$ (dashed). The inner of the two curves has only three right angles, contradicting Lemma~\ref{lem:four-corners}.}
\label{fig:triangle}
\end{figure}

\begin{lemma}
\label{lem:3way-helly}
Let $A$, $B$, and $C$ be orthoconvex regions of a Manhattan orbifold and suppose that each of the three intersections $A\cap B$, $B\cap C$, and $A\cap C$ is nonempty. Then the intersection of all three regions $A\cap B\cap C$ must be nonempty.
\end{lemma}

\begin{proof}
Suppose for a contradiction that the three pairwise intersections are nonempty but that the three-way intersection is not. As in the proof of the previous lemma, we may assume that the boundary of $A$ partitions the boundaries of $B$ and $C$ into \emph{outer arcs} exterior to $A$ together with components interior to $A$. By Lemma~\ref{lem:ocintcon}, the intersection of $B$ and $C$ can be bounded only by one arc from $B$ and one arc from $C$; let these arcs be called $b$ and $c$ respectively.

There are two simple closed curves that follow part of arc $c$ from one of its two endpoints until the first point of $b\cap c$, and then continue around arc $b$ (following the boundary of $B\cup C$ near the intersection point) until its endpoint on the boundary of $A$, and that finally continue around the boundary of $A$ (following the boundary of $B\cup A$ near the endpoint of $b$ until returning to the initial endpoint of $c$ (Figure~\ref{fig:triangle}). At least one of these two curves has an interior that is disjoint from $A\cup B\cup C$; for, if instead $A$ is contained in the interior of one of these curves, then it contains as well the portion of the other curve that follows the boundary of $A$, and it is not possible for both curves to contain each other. But if one of these curves has an interior that is disjoint from $A\cup B\cup C$, then by an argument similar to the one in the proof of the previous lemma it forms an orthogonal polygon with at most three right angles at the three points where it changes from following one boundary curve to following another. This violation of Lemma~\ref{lem:four-corners} shows that the pairwise intersections cannot be disjoint, so the triple intersection must be nonempty.
\end{proof}

A family of sets is a \emph{Helly family} if any pairwise intersecting subfamily has a common intersection.

\begin{lemma}
\label{lem:oc-helly}
The family of orthoconvex regions of any Manhattan orbifold is a Helly family.
\end{lemma}

\begin{proof}
Let $F$ be a pairwise intersecting family of orthoconvex regions; we must show that $F$ has a common intersection. Let $R_0$ be any region in $F$; by Lemma~\ref{lem:finite-singular}, $R_0$ contains a finite number of singularities of the orbifold. After having chosen a region $R_i$, if $R_i$ contains no singularities, the result follows by the Helly property for diagonally-aligned rectangles in the Manhattan plane. Otherwise, let $s$ be a singularity in $R_i$. If all other regions in the family contain $s$, then $s$ is a point of common intersection; otherwise let $T$ be a region in $F$ that does not contain $s$ and let $R_{i+1}=R_{i}\cap T$.  By Lemma~\ref{lem:ocintcon}, $R_{i+1}$ is itself an orthoconvex region. By Lemma~\ref{lem:3way-helly}, $R_{i+1}$ has a nonempty intersection with each other region in $F$, and it does not contain $s$. Each step of this type eliminates at least one singularity, while preserving the property that $R_i$ together with the sets in $F$ forms a pairwise intersecting family of orthoconvex regions. We started with finitely many singularities in $R_0$, so after finitely many steps the process of defining sets $R_i$ described above must terminate with a common intersection point.
\end{proof}

\section{The shape of a ball}

In order to show that Manhattan orbifolds are injective, we need to show that their metric balls form a Helly family. A ball is a set $B_r(p)=\{q\mid d(p,q)\le r\}$; the Helly property of these sets will follow from Lemma~\ref{lem:oc-helly} and from the following result, which shows that balls are orthoconvex.

\begin{figure}[t]
\centering\includegraphics[height=1.5in]{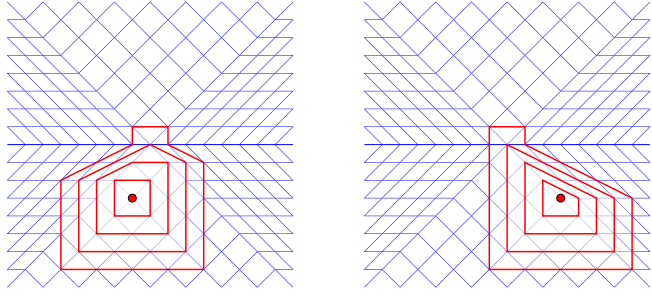}
\caption{A growing ball meets a cone point.}
\label{fig:ballxcone}
\end{figure}

\begin{lemma}
\label{lem:ball-oc}
Any ball in a Manhattan orbifold is orthoconvex.
\end{lemma}

\begin{proof}
To prove orthoconvexity of the ball $B_r(p)$, we consider the family of all balls $B_{r'}(p)$ for $r'\le r$. We show that there are finitely many combinatorially distinct shapes of balls in this family, and (by induction on the number of distinct shapes) that each is orthoconvex.

As a base case, when $p$ is not a singular point, all balls $B_\epsilon(p)$ for sufficiently small $\epsilon$ are isometric to a ball in the Manhattan plane, that is, a diagonally-aligned square. Such a shape is clearly orthoconvex. If $p$ is a cone point, the balls $B_\epsilon(p)$ for sufficiently small $\epsilon$ are isometric to the balls in a rectilinear cone, which take the form of a triangle in each quadrant of the cone together with a vertex with interior angle $\pi/2$ on each of the rays at which these quadrants are glued; again, such shapes are clearly orthoconvex.

We now consider how the shape of a ball can differ between $B_{r'}(p)$ and $B_{r'+\epsilon}(p)$, for sufficiently small values of $\epsilon$. The only configurations that can cause a change of shape in the ball between these two radii are those in which $B_{r'+\epsilon}(p)$ contains a singular point that $B_{r'}(p)$ does not, or those in which $B_{r'+\epsilon}(p)$ contains a point of the boundary of the orbifold that $B_{r'}(p)$ does not; in each case we can let $\epsilon$ be the minimum value possible that leads to this change. We note that it is not possible for the shape of the ball to change by an event in which two different parts of the boundary of the ball collide with each other, as the portion of boundary between the two colliding points would form an orthogonal polygon in which only one vertex has interior angle $\pi/2$, violating Lemma~\ref{lem:four-corners}.

We now describe in more detail each possible combinatorial change caused by the boundary of the ball reaching a singular point or a boundary point. We note that several such changes can happen at the same radius $r'+\epsilon$, but they can be considered independently of each other.

\begin{itemize}
\item If a ball meets a cone point of order $k$ along one of the edges of the correponding region (Figure~\ref{fig:ballxcone}, left) that cone point forms a vertex with interior angle $\pi$. For radii slightly larger than the radius at which the ball meets the cone point, there will be an additional $k-4$ vertices of interior angle $\pi/2$ near the cone point. For instance, Figure~\ref{fig:ballxcone}, left, shows a cone point of order six. The small nested balls have four vertices, the ball meeting the cone point has five vertices, and the largest ball has six vertices; note that some of the corners of the polygonal drawings of these balls are not vertices in the geometry of the depicted Manhattan orbifold.

\item If a vertex of a ball meets a cone point of order $k$ (Figure~\ref{fig:ballxcone}, right) no combinatorial change occurs until the radius grows larger than the radius at which the meeting occurs. For larger radii, there will be an additional $k-4$ vertices of interior angle $\pi/2$ near the cone point.

\begin{figure}[t]
\centering\includegraphics[height=1.5in]{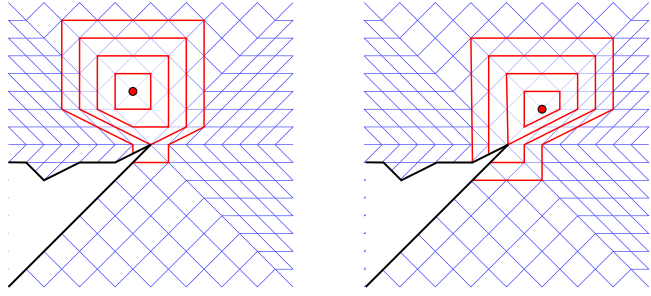}
\caption{A growing ball meets a cone inflection point forming a new intersection with the orbifold boundary.}
\label{fig:ballxinflex}
\end{figure}

\item If a ball meets an inflection point or a cone inflection point either along one of its edges or at a vertex that is not on the orbifold boundary for smaller radii, then for radii larger than the radius at which the ball meets this singularity, the ball's boundary includes a segment of the orbifold boundary (Figure~\ref{fig:ballxinflex}). There may also be additional vertices of interior angle $\pi/2$ near the singularity, depending on the relative angles of the ball boundary and orbifold boundary at the point singularity.

\item If a portion of the boundary of a ball, already including a segment of the orbifold boundary, meets a singularity on that segment of boundary (Figure~\ref{fig:ballxbdy}, left) then for larger radii there may be additional vertices of interior angle $\pi/2$ near the singularity, depending on the relative angles of the ball boundary and orbifold boundary at the singularity.

\item If a ball meets the orbifold boundary either along one of its edges or at a vertex that is not on the orbifold boundary for smaller radii, and does not meet an inflection point or cone inflection point  (Figure~\ref{fig:ballxbdy}, left), then the number of vertices of the ball interior to the orbifold is reduced by either one or two.
\end{itemize}

\begin{figure}[t]
\centering\includegraphics[height=1.5in]{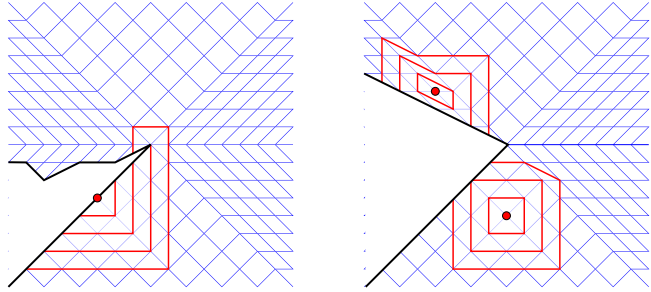}
\caption{A growing ball crosses a singularity on the boundary without meeting a new segment of boundary (left) or meets a new segment of boundary without crossing a singularity (right).}
\label{fig:ballxbdy}
\end{figure}

Thus, each possible combinatorial change between $B_{r'}(p)$ and $B_{r'+\epsilon}(p)$ preserves the orthoconvexity of the ball. Each change either includes within the ball an additional one of the finitely many singular points in $B_r(p)$, or reduces the finite number of vertices of the orthoconvex region, so only finitely many such changes are possible before a ball combinatorially equivalent to $B_r(p)$ will be reached. Therefore, $B_r(p)$ is orthoconvex.
\end{proof}

\section{Injectivity of Manhattan orbifolds}

We are now ready to prove our main result.

\begin{theorem}
Any Manhattan orbifold is injective.
\end{theorem}

\begin{proof}
By the result of \citet{AroPan-PJM-56} it suffices to prove that Manhattan orbifolds are path-geodesic and that their metric balls form a Helly family. The property of being path-geodesic was included in our definition of Manhattan orbifolds, and the Helly property of balls follows from Lemma~\ref{lem:oc-helly} and Lemma~\ref{lem:ball-oc}.
\end{proof}

\section{Tight spans of graphs}

As \citet{Isb-CMH-64} showed, any metric space $(X,d)$ can be isometrically embedded in a unique minimal injective space $T(X)$ called its \emph{injective envelope}, \emph{hyperconvex hull}, or \emph{tight span}. The tight span of a finite metric space can be defined as the set of functions $f:X\mapsto\R$, with the $L_\infty$ metric, satisfying the following properties:
\begin{itemize}
\item For every $p$ and $q$ in $X$, $f(p)+f(q)\ge d(p,q)$. In particular, taking $p=q$, $f(p)\ge 0$.
\item For every $p$ in $X$ there exists $q$ in $X$ such that $f(p)+f(q)=d(p,q)$.
\end{itemize}
This functions satisfying only the first of these two properties form a set $P(X)$ called the \emph{associated polytope} of $X$; the tight span can be viewed geometrically as the union of the bounded faces of this polytope.
Each point $p$ in $X$ corresponds to a function $f_p(q)=d(p,q)$, and it is straightforward to verify that this correspondence is an isometric embedding of $X$ into its tight span.
If $X$ is isometrically embedded into any injective space $S$, we can extend the embedding to $T(X)$: to find the point in $S$ corresponding to a function $f$, use the Helly property of balls to find an intersection point of the balls $B_{f(p)}(p)$.
On the other hand, $T(X)$ must contain an isometric copy of $T(Y)$ for every $Y\subset X$.

\begin{lemma}
Let $S$ be an injective metric space, let $X$ be a subspace of $S$, and let $X_i$ be a family of subsets of $X$. Suppose that for each $X_i$ there is a unique embedding of $T(X_i)$ in $S$ that is the identity on $X_i$, and suppose that the union of these copies of $T(X_i)$ covers $S$. Then $S$ is isometric to $T(X)$.
\end{lemma}

\begin{proof}
\label{lem:span-cover}
As discussed above, $T(X)$ must embed isometrically into $S$, as it does in every injective superset of $X$. But this embedding must include a copy of each $T(X_i)$, and therefore must cover all of $S$. Therefore $T(X)$ is isometric to $S$.
\end{proof}

Little is known about tight spans of graphs, but \citet{GooMou-DM-01} found the tight spans of certain graphs including cycles and hypercubes. Their results show that, for any $k>2$, the tight span of a $2k$-cycle is a hypercube of dimension $k$; in particular, no such graph can have a 2-manifold tight span.

\section{Squaregraphs}

\begin{figure}[t]
\centering\includegraphics[height=2in]{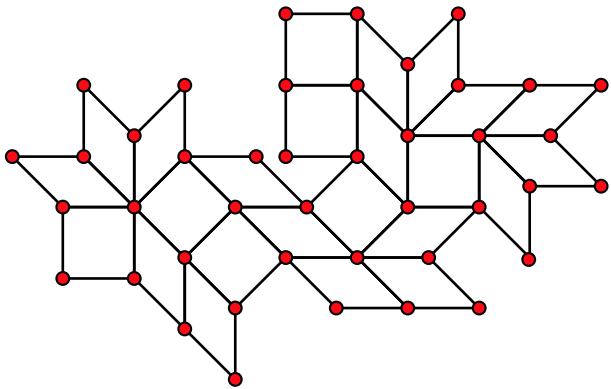}
\caption{A squaregraph.}
\label{fig:squaregraph}
\end{figure}

A \emph{squaregraph} \citep{CheDraVax-SODA-02,BanCheEpp-?} is a planar graph with an embedding in which all faces with the possible exception of the outer face are quadrilaterals and in which all vertices not part of the outer face have four or more incident edges. See Figure~\ref{fig:squaregraph} for an example.

To embed the squaregraph into an injective space, we associate each interior face of the squaregraph with a unit square of the Manhattan metric, and glue these squares together when the corresponding faces share an edge. The resulting space is clearly a Manhattan orbifold, hence injective. A squaregraph is a special case of a median graph~\citep{BanChe-05}, and this construction can be viewed as assigning the $L_1$ metric to the \emph{median polyhedral complex} \citep{BanChe-05,VdV-93} of the graph. It is known that the $L_\infty$ metric on the median polyhedral complex is injective~\citep{BanChe-05,MaiTan-PAMS-83,VdV-PAMS-98} however it does not contain an isometric copy of the original squaregraph, as vertices on opposite corners of a face are mapped to a unit distance apart. Using the $L_1$ metric for the median polyhedral complex avoids this problem.

\begin{lemma}
\label{lem:sqg-iso}
The map from each vertex of a squaregraph to the corresponding point of its $L_1$ median polyhedral complex is isometric.
\end{lemma}

\begin{proof}
Any path in the surface can be transformed, one square at a time starting with the vertex at the endpoint of the path, to a path of equal length that avoids the interior of any square. Therefore, the distance between vertices in the original graph equals the length of the shortest path in the surface.
\end{proof}

\begin{theorem}
\label{thm:squaregraph-span}
The $L_1$ median polyhedral complex of a squaregraph is isometric to its tight span.
\end{theorem}

\begin{proof}
We have already seen that this complex is injective, and that it contains an isometric copy of the graph. For each interior face of the graph, we associate a set $X_i$ consisting of the four vertices of that face; the tight span $T(X_i)$ embeds uniquely into the complex as the square corresponding to that face.
Therefore, by Lemma~\ref{lem:span-cover}, the complex is the tight span of the whole graph.
\end{proof}

It is tempting to try to extend this result to more general median graphs, however in general the tight span of a median graph cannot be formed from its median complex. In particular, as \citet{GooMou-DM-01} showed, the tight span of a cube graph is not homeomorphic to a geometric cube, but rather to a four-dimensional polytope.

\section{An infinite squaregraph}

\begin{figure}[t]
\centering\includegraphics[height=3.5in]{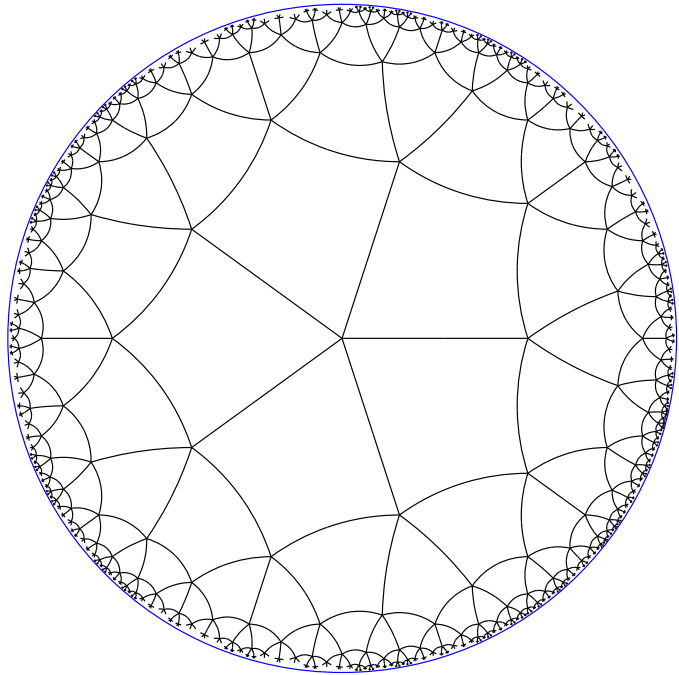}
\caption{The $\{4,5\}$ tesselation of the hyperbolic plane. Image produced using a Java applet written by Don Hatch, http://www.plunk.org/$\sim$hatch/HyperbolicApplet/.}
\label{fig:45tess}
\end{figure}

Figure~\ref{fig:45tess} shows the $\{4,5\}$ tesselation of the hyperbolic plane, a tiling of the plane by congruent squares meeting five at a corner. As is evident from the figure, this tesselation is an infinite squaregraph. Theorem~\ref{thm:squaregraph-span} applies equally well to infinite squaregraphs as it does to finite squaregraphs, so the tight span of this squaregraph is an unbounded Manhattan orbifold formed by replacing the hyperbolic metric within each square of the tiling by a unit square with the Manhattan metric.

With a suitable scaling factor, this surface can be mapped one-to-one onto the hyperbolic plane in such a way that the map distorts any distance by a small constant factor, analogously to the way the Manhattan metric on the plane distorts the Euclidean distance by at most a factor of $\sqrt 2$. It is possible that this injective approximation to the Euclidean metric may find some applications in geometric approximation algorithms for the hyperbolic plane \citep{Epp-TALG-?,KraLee-FOCS-06}, analogous to algorithms that approximate the Euclidean plane by a Manhattan metric. For instance, \citet{BeiChrLar-TCS-02} develop a 3-competitive 3-server algorithm in the Manhattan plane, from which it immediately follows that the same algorithm is $3\sqrt 2$-competitive in the Euclidean plane. If this algorithm could be generalized from the Manhattan plane to other Manhattan orbifolds, it could be made to apply in the same way to the hyperbolic plane.

\section{Beyond squaregraphs}

\begin{figure}[t]
\centering\includegraphics[height=1.5in]{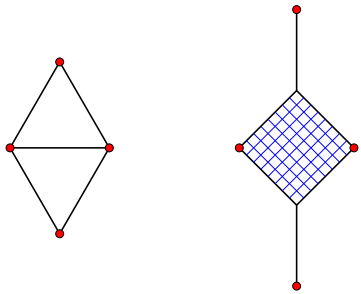}
\caption{A four-vertex graph with two triangles, and its tight span.}
\label{fig:two-triangles}
\end{figure}

Not every graph with a two-dimensional tight span is a squaregraph. For instance, if $G$ is a complete graph minus one edge (Figure~\ref{fig:two-triangles}) its tight span consists of a square with side length $1/2$ in the Manhattan plane, together with two one-dimensional \emph{whiskers} consisting of length-$1/2$ paths connected to two opposite corners of the square. The vertices of the graph embed in the tight span as the endpoints of the whiskers and the remaining two corners of the square.

We can use this square-and-whiskers construction to form the tight spans of a much larger class of graphs, including some highly nonplanar graphs.
A \emph{kinggraph} \citep{CheDraVax-SODA-02} is a graph formed from a squaregraph by adding edges connecting the diagonals of each of the squaregraph's quadrilateral faces. For instance, the graph formed in this way from an $8\times 8$ grid graph represents the possible moves of a king on a chessboard. More generally, if $G$ is any graph embedded in the plane in such a way that each interior face has four or more edges and each interior vertex has degree at least four, then we define a \emph{cliquegraph} (Figure~\ref{fig:cliquegraph}) to be the graph formed by adding edges connecting any two vertices belonging to the same face in $G$.

\begin{theorem}
\label{thm:cliquegraph}
The tight span of any cliquegraph is a Manhattan orbifold together with possibly some length-$1/2$ whiskers attached to the boundary of the orbifold.
\end{theorem}

\begin{proof}
We form a Manhattan orbifold by associating a square of side length $1/2$ in Manhattan geometry with each internal edge of the planar graph $G$ from which the cliquegraph was formed. We associate two opposite corners of the square with the two vertices at the endpoints of the edges, and the other two opposite corners with the two faces on opposite sides of the edge; this association gives us a gluing rule for connecting the squares into a single surface. The requirements that $G$ have four edges per interior vertex or face imply that this surface is a Manhattan orbifold. For any vertex that is not the endpoint of an internal edge of the graph, we add an whisker of length $1/2$ connecting that vertex to the point associated with the face to which it belongs. This construction is illustrated in Figure~\ref{fig:cliquegraph}.

\begin{figure}[t]
\centering\includegraphics[height=2.5in]{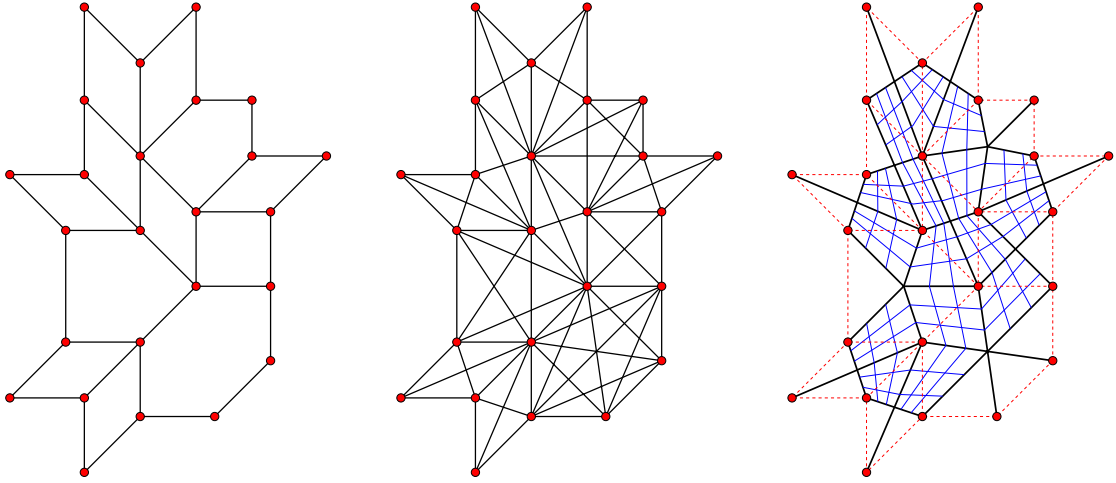}
\caption{A planar graph in which all interior vertices and faces have at least four edges, the corresponding cliquegraph, and its tight span.}
\label{fig:cliquegraph}
\end{figure}

It is clear that distances in the resulting surface are at most equal to distances in the graph.
By an argument similar to that of Lemma~\ref{lem:sqg-iso}, any path in the surface starting and ending at a vertex corresponds to a path with the same length and the same endpoints that avoids points in the interior of the glued-together squares; each adjacent pair of square edges in such a path has length one and connects two adjacent vertices of the cliquegraph, so distances in the surface equal distances in the graph. Thus, we have embedded the cliquegraph isometrically into an injective space.

Each whisker of the construction is part of the tight span of the whisker endpoint and another vertex on the same face, a tight span that embeds uniquely into our constructed surface. Each of the squares from which our surface is formed is part of the tight span of a two-triangle graph formed from the edge corresponding to the square together with two other vertices of the two faces on opposite sides of that edge; again, this tight span embeds uniquely into our surface. Therefore, by Lemma~\ref{lem:span-cover}, our construction is the tight span of the whole graph.
\end{proof}

\begin{figure}[t]
\centering\includegraphics[height=1.5in]{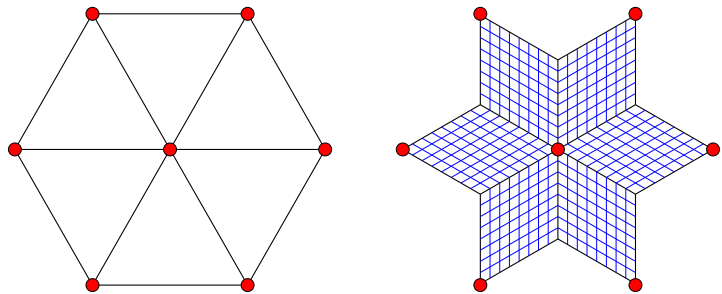}
\caption{A wheel and its tight span.}
\label{fig:wheelspan}
\end{figure}

A \emph{wheel} is a planar graph formed from a cycle and one additional vertex, called the hub of the wheel; it has edges connecting every vertex of the cycle to the hub (Figure~\ref{fig:wheelspan}).

\begin{theorem}
\label{thm:wheelspan}
If $G$ is a wheel with more than four vertices, then the tight span of $G$ is a Manhattan orbifold.\
\end{theorem}

\begin{proof}
For each vertex of the cycle in the wheel, we form a square in Manhattan geometry with side length $1/2$; we glue these squares together at a common cone point, placing the hub of the graph at this cone point and the cycle vertices of the graph on each square diagonally opposite the hub, as shown in Figure~\ref{fig:wheelspan}. As in Theorem~\ref{thm:cliquegraph}, each of the squares from which our surface is formed is part of the tight span of two adjacent triangles in the wheel, and each such tight span embeds uniquely into the overall surface, so by Lemma~\ref{lem:span-cover}, our construction is the tight span of the whole graph.
\end{proof}

A four-vertex wheel is just a clique $K_4$, the tight span of which consists of four length-$1/2$ whiskers connected at a common vertex. Thus it, too, can be isometrically embedded into a Manhattan orbifold, although its tight span is only one-dimensional.

With squaregraphs, cliquegraphs, and wheels, we have not exhausted the set of graphs that may be isometrically embedded into Manhattan orbifolds. A 5-vertex cycle, for instance, may be isometrically embedded into an order-5 rectilinear cone; its tight span is formed by five squares, with side length $1/2$ in Manhattan geometry, glued together at a common vertex. A \emph{house} formed by adding a single diagonal to the 5-cycle may be embedded directly into the Manhattan plane; its tight span is a unit Manhattan square together with an whisker attached to the midpoint of one of the square's edges. It would be of interest to characterize the isometric subgraphs of Manhattan orbifolds, but such a result is beyond the scope of the present work.

\section{Greedy embedding of graphs}

\begin{figure}[t]
\centering\includegraphics[height=2.5in]{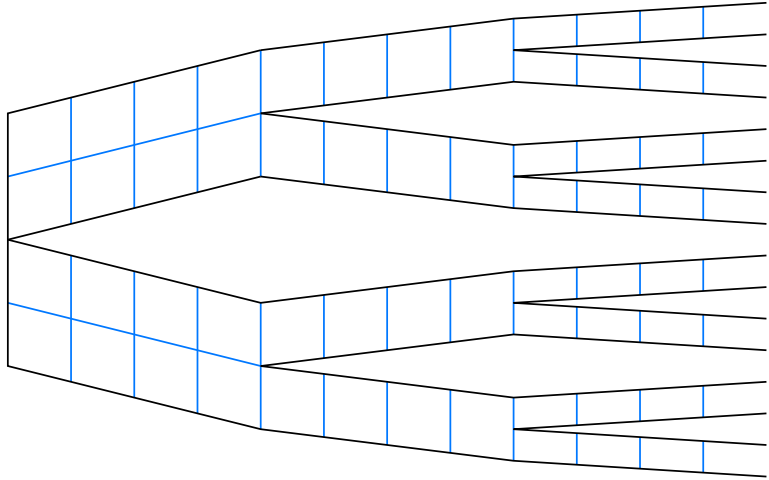}
\caption{The tight span of the dyadic tree metric space.}
\label{fig:dyadic}
\end{figure}

A \emph{greedy embedding} of a connected graph into a metric space is a function $f$ that maps the vertices of the graph to points of the metric space, in such a way that, for any two vertices $v\ne w$, there exists a neighbor $u$ of $v$ with $d(u,w)<d(v,w)$. Embeddings of this type may be used to define a simple greedy routing strategy for transmitting messages from one vertex to another in the graph: whenever a message with final destination $w$ reaches a vertex $v$, it is routed to a neighbor $u$ that takes it closer to its destination \citep{PapRat-TCS-05}. Greedy routings and greedy embeddings have been much studied in computer science; two of many results in this area is that every connected graph has a greedy embedding in the hyperbolic plane \citep{Kle-INFOCOM-07} and that descriptions of the vertex locations in such a representation can be encoded succinctly \citep{EppGoo-GD-08}. As part of their proof of the latter result, Eppstein and Goodrich described an embedding in a more abstract metric space that they called the \emph{dyadic tree metric space}. In this section we describe this space and show that its tight span is (essentially) a Manhattan orbifold. As a consequence of the results here and in \citet{EppGoo-GD-08}, it follows that every graph has a greedy embedding into a Manhattan orbifold.

The definition of the dyadic tree metric space is based on an embedding $f$ of an infinite binary tree $B$ into the unit interval: the root of the tree is mapped to the number $1/2$, the left subtree is mapped in the same fashion into the subinterval $[0,1/2]$, and the right subtree is mapped into the same fashion into the subinterval $[1/2,1]$. Thus, a vertex at level $i$ of tree $B$ is mapped into a dyadic rational number of the form $h/2^{i+1}$. A point in the dyadic tree metric space is defined to be a pair $(x,y)$, where $x$ and $y$ are two nodes in tree $B$ such that $x$ is an ancestor of $y$. The distance between points $(x,y)$ and $(x',y')$ in the dyadic tree metric space is defined to be $d(x,x')+|f(y)-f(y')|$, where distance in $B$ is an integer measuring the number of links in the shortest path between any two vertices.

For any fixed $x$ at level $i$ of tree $B$ and variable $y$, the points $(x,y)$ are isometric to the dyadic rationals (fractions with power-of-two-denominator) in the open interval $(f(x)-2^{-1-i},f(x)+2^{-1-i})$, so their tight span is this interval itself. We may form a path-geodesic metric space by connecting the interval formed in this way for a vertex $x$ to the interval for the parent of $x$ by an $L_1$-metric rectangle with side lengths $1$ and $2^{-i}$, with the intervals for $x$ and its parent identified with the short sides of this rectangle, as depicted in Figure~\ref{fig:dyadic}. These rectangles are the tight spans of the two intervals they connect. The resulting space satisfies the definition of a Manhattan orbifold except for having a single articulation point at the midpoint of the unit interval for the tree root. Thus, it is an injective space, the tight span of the dyadic tree metric space.

Since every finite undirected graph has a greedy embedding into the dyadic tree metric space \citep{EppGoo-GD-08}, it also has a greedy embedding into this tight span, and therefore into a Manhattan orbifold.

\raggedright
\bibliography{injective}

\end{document}